\documentclass[11pt]{article}       
\usepackage{graphicx}
\usepackage{mathptmx}      
\usepackage{times}

\usepackage{mathtools}
\usepackage{amsmath,amssymb,amsthm}
\usepackage{subfig}

\usepackage{bbm}
\usepackage{algorithmic}
\usepackage{algorithm}
\usepackage{color}
\usepackage[title]{appendix}

\newcommand{\Rb}{\mathbbm{R}}      

\newcommand{\Dc}{\mathcal{D}}
\newcommand{\Qc}{\mathcal{Q}}

\newcommand{\Eb}{\mathbbm{E}}
\newcommand{\Fc}{\mathcal{F}}

\newcommand{\Kc}{\mathcal{K}}

\newcommand{\Pc}{\mathcal{P}}

\newcommand{\Uc}{\mathcal{U}}

\newcommand{\Vc}{\mathcal{V}}
\newcommand{\Xc}{\mathcal{X}}

\newcommand{\1}{\mathbbm{1}}

\newcommand{\argmin}{\mathop{\rm argmin}}

\newcommand{\diag}{\mathop{\rm diag}}
\newcommand{\proj}{\text{\rm Proj}}
\newcommand{\conv}{\mathop{\rm conv}}
\newcommand{\range}{\mathop{\rm range}}

\newtheorem{assumption}{Assumption}

\newtheorem{theorem}{Theorem}
\newtheorem{lemma}[theorem]{Lemma}

\newenvironment{tightitemize}{%
    \list{{\textup{$\bullet$}}}{\settowidth\labelwidth{{\textup{\qquad}}}
    \leftmargin\labelwidth \advance\leftmargin\labelsep
    \parsep 0pt plus 1pt minus 1pt \topsep 3pt \itemsep 3pt
    }}{\endlist}

\newenvironment{tightenumerate}[1]{%
    \list{{\bfseries\textup{\arabic{enumi}.}}}{\settowidth\labelwidth{{\textup{(#1)}}}
    \leftmargin 0pt \advance\leftmargin\labelsep \itemindent \parindent
    \parsep 0pt plus 1pt minus 1pt \topsep 0pt \itemsep 0pt
    \usecounter{enumi}}}{\endlist}

   \makeatletter
\newcommand*\rel@kern[1]{\kern#1\dimexpr\macc@kerna}
\newcommand*\widebar[1]{%
  \begingroup
  \def\mathaccent##1##2{%
    \rel@kern{0.8}%
    \overline{\rel@kern{-0.8}\macc@nucleus\rel@kern{0.2}}%
    \rel@kern{-0.2}%
  }%
  \macc@depth\@ne
  \let\math@bgroup\@empty \let\math@egroup\macc@set@skewchar
  \mathsurround\z@ \frozen@everymath{\mathgroup\macc@group\relax}%
  \macc@set@skewchar\relax
  \let\mathaccentV\macc@nested@a
  \macc@nested@a\relax111{#1}%
  \endgroup
}
\makeatother

\begin{document}

\title{Risk-Averse Learning by Temporal Difference Methods}

\author{\"{U}mit K\"{o}se and Andrzej Ruszczy\'nski\footnote{Department of Management Science and Information Systems,
       Rutgers University,
       Piscataway, NJ 08854, USA; email: \{uek1,rusz\}@rutgers.edu; the order of authors is alphabetical}
       }

\date{March 1,2020}
\maketitle

\begin{abstract}
We consider reinforcement learning with performance evaluated by a dynamic risk measure.
We construct a projected risk-averse dynamic programming equation and study its properties. Then
we propose risk-averse counterparts of  the methods of temporal differences  and we prove their convergence with probability one.
We also perform an empirical study on a complex transportation problem.
\emph{AMS:}
  49L20, 62L20, 90C39
\end{abstract}

\section{Introduction}

The objective of this paper is to propose and analyze new risk-averse reinforcement learning methods for Markov Decision Processes (MDPs).
Our goal is to combine the efficacy of the methods of temporal differences with the robustness of dynamic risk measures,
and to provide a rigorous mathematical analysis of the methods.

MDPs are well-known models of stochastic sequential decision problems, covered in multiple monographs
\cite{Bellman1957a, Howard1960, Puterman1994,bertsekas2017dynamic}, and having countless applications.
In the classical setting, the goal of an MDP is to find a policy minimizing the expected cost over a finite or infinite horizon.
Traditional MDP models, although effective for small to medium size problems, suffer from the curse of dimensionality in problems with large state space. Approximate dynamic programming approaches try to tackle the curse of dimensionality and provide an approximate solution of an  MDP (see \cite{Powell2011} for an overview). Such methods usually involve value function approximations, where the value of a state of the Markov process is approximated by a simple, usually linear, function of some selected features of the state \cite{Bellman1963}.

Reinforcement learning methods \cite{Sutton1998, Powell2011} involve simulation or observation of a Markov process to approximate the value function and learn the corresponding policies. The first studies attempted to emulate neural networks and biological learning processes, learning by trial and error \cite{Minsky1954, Farley1954}. Some learning algorithms, such as Q-Learning \cite{Watkins1989, Watkins1992} and SARSA \cite{Rummery1994}, follow this idea. One of core approaches in reinforcement learning is the method
of temporal differences \cite{Sutton1988}, known as TD($\lambda$). It uses differences between the values of the approximate value function at successive states to improve the approximation, concurrently with the evolution of the system. TD($\lambda$) is a continuum of algorithms depending on a parameter $\lambda \in [0,1]$ which is used to exponentially weight past observations.  Consequently, related methods such as Q($\lambda$) \cite{Watkins1989, Peng1993, Peng1994, Rummery1995} and SARSA($\lambda$) were developed \cite{Rummery1994, Rummery1995}. The methods of temporal differences have been proven to converge in the mean in \cite{Dayan1992} and almost surely by several studies,
with different degrees of generality and precision
\cite{Peng1993,Dayan1994,Tsitsiklis1994, Jaakkola1994,tsitsiklis1997analysis}.

We introduce risk models into temporal difference learning.
In the extant literature, three basic approaches to risk aversion in  MDPs have been employed:
utility functions (see, e.g., \cite{Jaquette1973,Jaquette1975,Chung1987,Denardo1979,Fleming1999,bauerle2013more,jaskiewicz2013persistently}), mean--variance models (see. e.g.,
\cite{White1988,Filar1989,Mannor2013,arlotto2014markov,ChenLiZhao2014}), and entropic (exponential) models (see, e.g.,
\cite{Howard1971,Marcus1997,bielecki1999risk,Coraluppi1999,DiMasi1999,Levitt2001,bauerle2013more}). Our research
is rooted in the theory of dynamic measures of risk, which has been intensively developed in the last 15 years (see
 \cite{Scandolo:2003,Riedel:2004,roorda2005coherent,follmer2006convex,CDK:2006,RuszczynskiShapiro2006b,ADEHK:2007,PflRom:07,KloSch:2008,jobert2008valuations,cheridito2011composition}
 and the references therein).

In \cite{Ruszczynski2010Markov},  we introduced the class of Markov dynamic risk measures, specially tailored for the MDPs.
It allowed for the development of dynamic programming equations and corresponding solution methods, generalizing the well-known results for the expected value problems. Our ideas were successfully extended to undiscounted problems in  \cite{CavusRuszczynski2014a,CavusRuszczynski2014b},
partially observable and history-dependent systems in \cite{fan2018risk,fan2018process}, and further generalized in \cite{lin2013dynamic,shen2013risk}.

A number of works introduce models of risk into reinforcement learning: exponential utility functions \cite{Borkar2001, Borkar2002} and  mean-variance models \cite{Tamar2012, Prashanth2014}. Few later studies propose heuristic approaches involving coherent risk measures and their mean-risk counterparts \cite{Chow2014, Tamar2017}; these studies employ policy gradients and use them in actor-critic type algorithms.
Distributed policy gradient methods with risk measures were proposed in \cite{ma2017risk}.
 Model-related uncertainties are discussed in \cite{Tamar2014}.

In this paper, we use Markov risk measures of \cite{Ruszczynski2010Markov} in conjunction with linear approximations of the value function.  Our contributions can be summarized as follows:
\begin{tightitemize}
\item A projected risk-averse dynamic programming equation and analysis of its properties (\S \ref{sec:projec});
\item A risk-averse method of temporal differences (\S \ref{sec:methods})
and proof of its convergence (\S \ref{s:TD0-convergence});
\item A multistep risk-averse method of temporal differences (\S  \ref{s:TDL}) and its convergence proof (\S \ref{s:TDL-convergence});
\item An empirical study comparing the efficacy of the methods (\S \ref{s:empirical}).
\end{tightitemize}

\section{The Projected Risk-Averse Dynamic Programming Equation}
\label{sec:projec}
We consider a Markov decision process (MDP) with a finite state space $\mathcal{X}= \{1,\dots,n\}$,  finite action sets $\mathcal{U}(i)$ for all $i \in \mathcal{X}$, controlled transition probabilities $P_{ij}(u)$ where $i,j \in \mathcal{X}$ and $u \in \mathcal{U}(i)$, and one-step cost function $c(i,u)$, where $i \in \mathcal{X}$ and $u \in \mathcal{U}(i)$. For a discount factor $\alpha\in (0,1)$ and any non-anticipative policy $\pi$ for determining controls $u_t\in \Uc(i_t)$, $t=0,1,2,\dots$,  the expected discounted cost
\[
v^\pi(i) = \Eb\Big[ \sum_{t=0}^\infty \alpha^t c(i_t,u_t)\,\Big| \,i_0=i\Big],
\]
is finite.
For every Markovian policy $\pi$, the value function associated with this policy
satisfies the linear equation
\[
v^\pi(i) =  c(i,\pi(i)) + \alpha \sum_{j\in \Xc}P_{ij}(\pi(i)) \,v^\pi(j), \quad  i \in \mathcal{X}.
\]
Viewing $v^\pi$ as a vector, and defining the vector $c^\pi$ with elements $c(i,\pi(i))$, $i\in \Xc$, and the matrix $P^\pi$ with entries $P_{ij}(\pi(i))$, $i,j\in \Xc$, we can compactly
write the policy evaluation equation as
\begin{equation}
\label{policy-neutral}
v^\pi = c^\pi + \alpha P^\pi v^\pi.
\end{equation}
In \cite{Ruszczynski2010Markov}, in a more general setting in a Polish space $\Xc$,  Markov risk measures for cost evaluation in an MDP were introduced. In a finite-horizon setting,
a Markov risk measure evaluates the sequence of discounted costs $\alpha^t c(x_t,u_t)$, $t=0,1,2,\dots,T$, under a Markov policy $\pi$, in a recursive way. Denoting by $\rho^\pi_{t,T}(i)$ the risk
of the system starting from state $i$ at time $t$, we have
\begin{equation}
\label{DP-risk-finite}
\rho^\pi_{t,T}(i) = c^\pi_i + \alpha \sigma_i\big(P^\pi_i, \rho^\pi_{t+1,T}(\cdot)\big),\quad i\in \Xc, \quad t=0,1,\dots,T-1,
\end{equation}
with $\rho^\pi_{T,T}(i) = c^\pi_i$, $i\in \Xc$. In equation \eqref{DP-risk-finite}, the operator $\sigma:\Xc \times  \Pc(\Xc)\times \Vc\to \Rb$,
where $\Pc(\Xc)$ is the space of probability measures on $\Xc$ and $\Vc$ is the space of bounded functions on $\Xc$, is a \emph{transition risk mapping}. It
can be interpreted as risk-averse analog of the conditional expectation. Its first argument is the state $i$ (which we write as a subscript). The second argument, the vector $P^\pi_i$, is the $i$th row of the matrix $P^\pi$: the probability distribution of the state following $i$ under the policy $\pi$. The last argument, the function $\rho^\pi_{t+1,T}(\cdot)$, is the risk of running the system from the next state
in the time interval from $t+1$ to $T$.
The transition risk mapping is a special case of a \emph{risk form}: a generalization of a risk measure introduced in \cite{dentcheva2019risk} to accommodate the dependence of measures of risk on the underlying
probability distribution. In the case of controlled Markov systems, this dependence is germane for the analysis.

As in \cite{Ruszczynski2010Markov}, we assume that for each $i\in \Xc$ and each $P_i^\pi\in \Pc(\Xc)$, the transition risk mapping $\sigma_i(p,\cdot)$, understood as a function
its last argument, satisfies the axioms of a coherent measure of risk \cite{ArtznerDelbaenEberEtAl1999}. In the axioms below we suppress the argument $P^\pi_i$, focusing
on the dependence on the third argument, a function of a state:
\begin{tightitemize}
\item[\textbf{Convexity:}] $\sigma_i(\alpha v + (1-\alpha) w) \leq \alpha \sigma_i(v) + (1-\alpha) \sigma_i(w)$,  $\forall\, \alpha \in [0,1]$, $\forall\, v,w\in \Vc$;
\item[\textbf{Monotonicity:}] If $v \leq w$ (componentwise) then $\sigma_i(v) \leq \sigma_i(w)$;
\item[\textbf{Translation equivariance:}] $\sigma_i(v + \beta\1) = \sigma_i(v) + \beta$, for all $\beta \in \Rb$;
\item[\textbf{Positive homogeneity:}] $\sigma_i(\beta v) = \beta \sigma_i(v)$, for all $\beta \geq 0$.
\end{tightitemize}
Under these conditions, one can pass to the limit with $T\to \infty$ in \eqref{DP-risk-finite} and prove the existence of an {infinite-horizon discounted risk measure} \cite{Ruszczynski2010Markov}
\[
\rho^\pi_{0,\infty}(i) = \lim_{T\to\infty}\rho^\pi_{0,T}(i) , \quad i\in \Xc.
\]
We still denote its value at state $i$ by $v^\pi(i)$; it will never lead to misunderstanding. The policy value $v^\pi(\cdot)$ satisfies the {risk-averse policy evaluation equation}:
\[
v^\pi(i) = c^\pi(i) + \alpha \sigma_i\big(P^\pi_i, v^\pi(\cdot)\big), \quad i\in \Xc.
\]
We introduce the space $\Qc$ of transition kernels on $\Xc$, define a vector-valued \emph{transition risk operator} $\sigma: \Qc \times \Vc \to \Vc$, with components $\sigma_i(P^\pi_i,\cdot)$, $i\in \Xc$, and rewrite the last equation
in a way similar to \eqref{policy-neutral}:
\begin{equation}
\label{policy-risk}
v^\pi  = c^\pi + \alpha \sigma (P^\pi,v^\pi).
\end{equation}
The only difference between \eqref{policy-neutral} and \eqref{policy-risk} is that the matrix $P^\pi$ has been replaced by a convex operator $\sigma$ (which still depends on $P^\pi$).
The risk-neutral case is a special case of \eqref{policy-risk} with $\sigma (P,v)=P v$.
References \cite{Ruszczynski2010Markov,CavusRuszczynski2014a,CavusRuszczynski2014b,fan2018process} outline the theory,  provide examples   and applications.

Coherent risk measures admit a dual representation \cite{RuszczynskiShapiro2006a}, which in our case can be stated as follows. For every $i\in \Xc$ a convex, closed and bounded set $A_i(P_i^\pi)$
of probability measures on $\Xc$ exists, such that
\begin{equation}
\label{sigma-dual}
\sigma_i(P^\pi_i,v) = \max_{\mu\in A_i(P_i^\pi)} \langle \mu, v \rangle,\quad v\in \Vc.
\end{equation}
In a risk-neutral case, the set $A_i(P_i^\pi) = \partial \sigma_i(P^\pi_i,0) $ contains only one element, $P_i^\pi$, but in general it is larger and has $P_i^\pi$ as one of its elements, provided we always have $\sigma_i(P^\pi_i,v) \ge P^\pi_i v$. The multifunction $A:\Xc\to\Pc(\Xc)\rightrightarrows\Pc(\Xc)$ is called the \emph{risk multikernel}. Every $\mu\in A_i(P_i^\pi)$ is absolutely continuous with respect to $P_i^\pi$.


While equation \eqref{policy-risk} can be solved by a nonsmooth Newton's method and the resulting evaluation used in a policy iteration method \cite{Ruszczynski2010Markov},
all these techniques become impractical, when the size of the state space is very large.

An established approach to such a situation in expected value models is to assume that each state $i\in \Xc$ has a number of relevant \emph{features} $\varphi_j(i)\in \Rb$, $j=1,\dots,m$, where $m \ll n$, and that the value $v^\pi(i)$ of a state can be approximated by a linear combination of its features:
\begin{equation}
\label{v-tilde-linear}
v^\pi(i) \approx \widetilde{v}^\pi(i)  = \sum_{j = 1}^m r_j\varphi_j(i), \quad i \in \mathcal{X}.
\end{equation}
From now on, we suppress the superscript $\pi$, because most of our considerations focus on evaluating a fixed policy.
We define the matrix of the features of all states, namely
\begin{align*}
\varPhi = \begin{bmatrix}\varphi^{\top}\!(1)\\ \varphi^{\top}\!(2)\\ \vdots\\ \varphi^{\top}\!(n)\end{bmatrix}.
\end{align*}
 Now we can write our approximation as $v \approx \tilde{v} = \varPhi r$. Similar to the expected value case, if we attempt to emulate \eqref{policy-risk} with the approximate value function, we may observe that the right hand side of the equation, $c + \alpha \sigma(P,\varPhi r)$, may not be represented as a linear combination of the features. Therefore, we need to project this vector on the subspace spanned by the features, $\range(\varPhi)$. Accordingly, we define a projection operator, $L: \mathcal{V} \rightarrow \range(\varPhi)$, and formulate the projected risk-averse approximate dynamic programming equation:
\begin{align}
\varPhi r  = L\big(c + \alpha \sigma(P,\varPhi r)\big). \label{projeq}
\end{align}
Still following the expected value case, we assume that the Markov system under policy $\pi$ is ergodic, and we denote its vector of stationary probabilities by $q$. We define
the projection operator using the following scalar product and the associated norm:
$\langle v,w\rangle_q = \sum_{i=1}^n q_i v_iw_i$, $\|w\|_q^2 = \langle w,w\rangle_q$.
Then
\begin{align}
L(w) = \argmin_{z \in \range(\varPhi)} ||z - w||_q, \quad w \in \mathcal{V}. \label{ldef}
\end{align}
The fundamental question is the existence and uniqueness of a solution of equation \eqref{projeq}.
This can be answered by establishing the contraction mapping property of the right hand side of \eqref{projeq}:
\begin{equation}
\label{D-operator}
\Dc(v) = L\big(c + \alpha \sigma(P,v)\big), \quad v\in \Vc,
\end{equation}
which would imply the existence and uniqueness of a solution of the equation
\begin{equation}
\label{DP-projected}
v = \Dc v.
\end{equation}
Crucial in this context  is the \emph{distortion coefficient} of the risk multikernel $A$:
\[
\varkappa =  \max \bigg\{ \frac{|\mu_{ij}-p_{ij}|}{p_{ij}} : \mu_i \in A_i(P_i^\pi),\  p_{ij} > 0,\  i,j \in \mathcal{X} \bigg\}.
\]
By definition, $\varkappa \ge 0$, with the value 0 corresponding to the risk-neutral model. We also recall
that for $p_{ij}=0$ we always have $m_{ij}=0$, for all $m_i\in A_i(P_i^\pi)$.
\begin{lemma}
\label{l:sigma-Lipschitz}
The transition risk operator satisfies for all $w,v\in \Vc$ the inequalities:
\begin{gather}
\label{sigma-L}
\| \sigma(P,w) - \sigma(P,v)\|_q \le \sqrt{1+\varkappa} \,\| w-v\|_q,
\intertext{and}
\label{sigma-delta}
\| \sigma(P,w) - \sigma(P,v) - P(w-v)\|_q \le \varkappa \,\| w-v\|_q .
\end{gather}
\end{lemma}
\begin{proof}
For brevity, we omit the argument $P$ of $\sigma(P,\cdot)$, because it is fixed. For every $i=1,\dots,n$,
by the mean value theorem for convex functions \cite{wegge1974mean,hiriart1980mean}, a point $\bar{v}^{(i)} = (1-\theta_i) v + \theta_i w$  exists, with
$\theta_i\in [0,1]$, and a subgradient $m_i\in \partial \rho_i(\bar{v}^{(i)})$ exists, such that
\[
\sigma_i(w) - \sigma_i(v) = \langle m_i, w -v \rangle.
\]
Since the subdifferential $\partial \rho_i(\cdot) \subseteq A_i$, we have $m_i \in A_i$.
Therefore, for a matrix $M$ having $m_i$, $i=1,\dots,n$, as its rows,
\begin{equation}
\label{mean-value}
\sigma(w) - \sigma(v) = M( w -v).
\end{equation}
As each $m_i$ is a probability vector,  Jensen's inequality  with $h=w-v$, and the equation $q^\top P = q^\top$ yield
\begin{multline}
\label{M-nonexpansive}
\|Mh\|_q^2  = \sum_{i\in \Xc} q_i \Big(\sum_{j\in \Xc} m_{ij}h_j\Big)^2 \le \sum_{i\in \Xc} q_i \sum_{j\in \Xc} m_{ij} h_j^2 \\
\le( 1+ \varkappa) \sum_{i\in \Xc} q_i \sum_{j\in \Xc} p_{ij} h_j^2 = (1+\varkappa) \sum_{j\in \Xc}q_j h_j^2  =(1+\varkappa) \|h\|_q^2.
\end{multline}
The last two relations imply \eqref{sigma-L}.
In a similar way, it follows from \eqref{mean-value} that
\begin{multline*}
\big\| \sigma(P,w) - \sigma(P,v) - P(w-v)\big\|^2_q  = \| (M-P)h\|^2_q
\le \sum_{i\in \Xc} q_i \Big(\sum_{j\in \Xc} |m_{ij}-p_{ij}| |h_j|\Big)^2\\
\le \varkappa^2 \sum_{i\in \Xc} q_i \Big(\sum_{j\in \Xc} p_{ij} |h_j|\Big)^2
\le \varkappa^2 \sum_{i\in \Xc} q_i \sum_{j\in \Xc} p_{ij} |h_j|^2
= \varkappa^2 \,\| w-v\|^2_q,
\end{multline*}
which is \eqref{sigma-delta}.
\end{proof}

We can now prove the existence and uniqueness of the solution of the risk-averse equation \eqref{DP-projected}.

\begin{theorem}
\label{contraction}
If $\alpha \sqrt{1+\varkappa} < 1$ then the equation \eqref{DP-projected} has a unique solution $v^*$.
\end{theorem}
\begin{proof}
We verify that the operator \eqref{D-operator} is a contraction mapping in the norm $\|\cdot\|_q$. The orthogonal projection $L$ is nonexpansive.
The operator $P$ is nonexpansive in the norm $\|\cdot\|_q$ as well (this is a special case of \eqref{M-nonexpansive} with $M=P$ and $\varkappa=0$).
The transition risk operator $\sigma(\cdot)$ multiplied by $\alpha$ is a con\-tract\-ion by Lemma \ref{l:sigma-Lipschitz}.
The assertion follows now from the Banach contraction mapping theorem.
\end{proof}
If $\varPhi$ has full column rank, equation \eqref{projeq} has a unique fixed point as well.

\section{The Risk-Averse Method of Temporal Differences}
\label{sec:methods}

We propose to solve \eqref{projeq} by a risk-averse analog of the classical method of temporal differences \cite{Sutton1988}.
We define ${v}^*$ to be the solution of equation \eqref{DP-projected} (which exists and is unique, if $\alpha\sqrt{1+\varkappa}<1$).

Consider the evolution of the
system under policy $\pi$, resulting in a random trajectory of states $i_t$, $t=0,1,2\dots$. At each time $t$, we have an approximation $r_t$ of a solution  of the equation
\eqref{projeq}. Let $\Fc_t$ be the $\sigma$-algebra defined by all observations gathered up to time $t$.

The difference between the left and the right hand sides of equation \eqref{projeq} with coefficient values $r_t$ and state $i_t$
is the \emph{risk-averse temporal difference}:
\begin{equation}
\label{hat-dt}
d_t = \varphi^{\top}\!(i_t)r_t - c(i_t) - \alpha \sigma_{i_t} (P_{i_t},\varPhi r_t), \quad t=0,1,2,\dots.
\end{equation}
Evidently, it cannot be easily computed or observed; this would require the evaluation of the risk $\sigma_{i_t} (P_{i_t},v)$ and thus consideration of \emph{all}
possible transitions from state $i_t$. Instead, we assume that we can observe a random estimate $\widetilde{\sigma}_{i_t}(P_{i_t},\cdot)$, such that
\begin{equation}
\label{sigma-tilde}
\widetilde{\sigma}_{i_t}(P_{i_t},\varPhi r_t) = \sigma_{i_t} (P_{i_t},\varPhi r_t) + \xi_{t},\quad t=0,1,2,\dots,
\end{equation}
with some random errors $\xi_{t}$. The conditions on $\{\xi_t\}$ will be specified later.
This allows us to define the observed risk-averse temporal differences,
\begin{equation}
\label{TD0-1}
\widetilde{d}_t = \varphi^{\top}\!(i_t)r_t - c(i_t) - \alpha \widetilde{\sigma}_{i_t} (P_{i_t},\varPhi r_t) , \quad t=0,1,2,\dots,
\end{equation}
and to construct the \emph{risk-averse temporal difference  method} as follows:
\begin{equation}
\label{TD0-2}
r_{t+1} =  r_t - \gamma_t \varphi(i_t)\, \widetilde{d}_t, \quad t=0,1,2,\dots.
\end{equation}

Before proceeding to the detailed convergence proof in the stochastic case, we analyze a deterministic model of the method, in which
the errors $\xi_t$ are ignored and the updates of the sequence $\{r_t\}$ are averaged over all states (with the distribution $q$).
We define the operator:
\begin{equation}
\label{U-operator}
U(r) = \Eb_{i \sim q}\big[ \varphi(i) \big(\varphi^{\top}\!(i)r - c(i)- \alpha \sigma_i(P_i,\varPhi r) \big)\big] =  \varPhi^{\top} Q \big[ \varPhi r -c - \alpha \,\sigma(P,\varPhi r  )\big].
\end{equation}
The deterministic analog of \eqref{TD0-1}--\eqref{TD0-2} reads:
\begin{equation}
\label{TD0-model}
\bar{r}_{t+1} = \bar{r}_t - \gamma\, U(r_t), \quad t=0,1,2,\dots,\quad \gamma>0.
\end{equation}
By the definition of the projection operator $L$, a point $r^*$ is a solution  of  \eqref{projeq}  if and only if
\[
r^* = \argmin_r\frac{1}{2}\big\|\varPhi r - \big(c+\alpha \sigma(P,\varPhi r^*)\big)\big\|_q^2.
\]
This occurs if and only if $r^*$ is a zero of $U(\cdot)$ and thus supports our idea of using the method \eqref{TD0-1}--\eqref{TD0-2}.
\begin{theorem}
\label{t:RTD-convergence}
If $\alpha\sqrt{1+\varkappa}<1$, then ${\gamma}_0>0$ exists, such that for all ${\gamma}\in (0,\gamma_0)$
the algorithm \eqref{TD0-model} generates a sequence $\{\bar{r}_t\}$ convergent to a point $r^*$ such that ${U}(r^*)=0$.
\end{theorem}
\begin{proof}
We shall show that for sufficiently small $\gamma>0$ the operator $I-\gamma U$ is a contraction.
For arbitrary $r'$ and $r''$, we have
\begin{multline*}
\big\| (r' - \gamma U(r'))- (r'' - \gamma U(r''))\big\|^2
=  \|r'-r''\|^2 \\
{} -2 \gamma \langle  r'-r'', \varPhi^{\top} Q \varPhi (r'-r'') \rangle
 + 2\gamma \alpha \big\langle  r'-r'', \varPhi^{\top} Q \big[ \sigma(P,\varPhi r')-\sigma(P,\varPhi r'')\big]\big\rangle\\
{ } + \gamma^2 \Big\| \varPhi^{\top} Q \varPhi (r'-r'')  - \alpha\varPhi^{\top} Q \big[ \sigma(P,\varPhi r')-\sigma(P,\varPhi r'')\big]\Big\|^2.
\end{multline*}
The last term (with $\gamma^2$) can be bounded by $\gamma^2 C \|\varPhi (r'-r'')\|_q^2$ where $C$ is some constant. Then
\begin{multline*}
\big\| (r' - \gamma U(r'))- (r'' - \gamma U(r''))\big\|^2
 \le  \|r'-r''\|^2 -2 \gamma \| \varPhi (r'-r'') \|_q^2 \\
{\qquad} + 2\gamma \alpha \big\langle \varPhi (r'-r''),  \sigma(\varPhi r')-\sigma(\varPhi r'')\big\rangle_q + \gamma^2 C \|\varPhi (r'-r'')\|_q^2.
\end{multline*}
The scalar product can be bounded by  \eqref{sigma-L}, and thus
\begin{align*}
\lefteqn{\big\| (r' - \gamma U(r'))- (r'' - \gamma U(r''))\big\|^2 }\quad\\
&  \le  \|r'-r''\|^2 -2 \gamma \| \varPhi (r'-r'') \|_q^2
 + 2\gamma \alpha \sqrt{1+\varkappa} \| \varPhi (r'-r'')\|^2_q + \gamma^2 C \|\varPhi (r'-r'')\|_q^2\\
&=  \|r'-r''\|^2 -2 \gamma  \bigg(1 - \alpha \sqrt{1+\varkappa} + \frac{\gamma C}{2} \bigg)\|\varPhi (r'-r'')\|_q^2.
\end{align*}
Since $\alpha \sqrt{1+\varkappa}< 1$, then using $0< \gamma< 2(1- \alpha\sqrt{1+\varkappa})/C $, we have
\begin{equation}
\label{semi-contraction}
\big\| (r' - \gamma U(r'))- (r'' - \gamma U(r''))\big\|^2  \le \|r'-r''\|^2 - \gamma \beta \|\varPhi (r'-r'')\|_q^2,
\end{equation}
with some $\beta>0$. In particular, setting $r' = \bar{r}_t$ and $r''=r^*$ for a solution $r^*$ of \eqref{projeq}, we obtain
the following relation between the successive iterates of the method \eqref{TD0-model}:
\begin{equation}
\label{descent-det}
\| \bar{r}_{t+1}-r^*\|^2 \le \| \bar{r}_{t}-r^*\|^2 - \gamma \beta  \| \varPhi(\bar{r}_{t}-r^*)\|_q^2.
\end{equation}
This immediately proves that the sequence $\{\bar{r}_t\}$ is bounded and $\varPhi \bar{r}_t \to \varPhi r^*$. Every accumulation point $\hat{r}$ of $\{\bar{r}_t\}$ must be then a solution
of equation \eqref{projeq}. Substituting this accumulation point for $r^*$ in the last inequality, we conclude that $\bar{r}_t \to \hat{r}$.
\end{proof}

If $\varPhi$ has full column rank, the solution $r^*$ is unique, because substituting another solution for $\bar{r}_t$ in \eqref{TD0-model} we obtain $r_{t+1} = r_t$, which
leads to a contradiction in \eqref{descent-det}.

\section{Convergence of the Risk-Averse Method of Temporal Differences}
\label{s:TD0-convergence}

We shall use the following result on convergence of deterministic nonmonotonic algorithms \cite{Nurminski1972}.
\begin{theorem}
\label{t:nurminski}
Let $Y^* \subset \Rb^m$. Suppose $\{r_t\} \subset \Rb^m$ is a bounded sequence  which satisfies the following assumptions:
\begin{tightitemize}
\item[A)] If a subsequence $\{r_t\}_{t\in \Kc}$ converges to $r'\in Y^*$, then
$\|r_{t+1}-r_t\|\to 0$,  as $t\to \infty$ , $t\in \Kc$;
\item[B)] If a subsequence $\{r_t\}_{t\in \Kc}$ converged to $r'\notin Y^*$, then $\varepsilon_0>0$ would exist such that
for all $\varepsilon\in(0,\varepsilon_0]$ and for all $k\in \Kc$, the index
$s(t,\varepsilon) = \min\big\{\ell \ge k: \|r_\ell-r_t\|>\varepsilon\big\}$
would be finite;
\item[C)] A continuous function $W:\Rb^m\to\Rb$ exists such that if $\{r_t\}_{t\in \Kc}$ converged to $r'\notin Y^*$
then $\varepsilon_1>0$ would exist such that for all $\varepsilon\in(0,\varepsilon_1]$ we would have
\[
\limsup_{t\in \Kc} W(r_{s(t,\varepsilon)})< W(r'),
\]
where $s(t,\varepsilon)$ is defined in B);
\item[D)] The set $\{W(r): r\in Y^*\}$ does not contain any segment of nonzero length.
\end{tightitemize}
Then the sequence $\{W(r_t)\}$ is convergent and all limit points of the sequence $\{r_t\}$ belong to $Y^*$.
\end{theorem}

We define the set of solutions of equation \eqref{projeq}:
\[
Y^* = \{r \in \Rb^m: \varPhi r = v^*\},
\]
where $v^*$ is the unique solution of \eqref{DP-projected}, provided $\alpha \sqrt{1+\varkappa}<1$
We shall show that the method \eqref{TD0-2} converges to $Y^*$, under the above-mentioned condition  and
some additional conditions on the stepsizes $\{\gamma_t\}$ and errors $\{\xi_t\}$.

We define $\Fc_t$ to be the $\sigma$-algebra generated by $\{i_0,r_0,\dots,i_t,r_t\}$, $t=0,1,\dots$, and
make the following assumptions about the stepsize and error sequences. We allow the stepsizes to be random.
\begin{assumption}
\label{a:gamma}
The sequence $\{\gamma_t\}$ is adapted to the filtration $\{\Fc_t\}$ and such that
\begin{tightenumerate}{iii}
\item[\rm(i)]  $\gamma_t >0$, $t=0,1,\dots$, and $\lim_{t\to\infty} \gamma_t = 0\quad \text{a.s.}$;
\item[\rm(ii)] ${\sum_{t=0}^\infty \gamma_t = \infty} \quad \text{a.s.}$;
\item[\rm(iii)] ${ \Eb \sum_{t=0}^\infty \gamma_t^2 < \infty}$;
\item[\rm(iv)] For any $\varepsilon>0$,
  $\displaystyle{\lim_{\raisebox{-0.75ex}{${\scriptstyle t_0\to \infty}$}\ }\sup_{\{T: \sum_{t=t_0}^T \gamma_t \le \varepsilon\}}}\
  \sum_{t=t_0}^T | \gamma_t-\gamma_{t+1}| = 0 \quad \text{a.s.}
$

 \end{tightenumerate} \vspace{0.5ex}
\end{assumption}
\begin{assumption}
\label{a:xi}
The sequence of errors $\{\xi_t\}_{t\ge 1}$ satisfies for $t=0,1,2\dots$ the conditions
\begin{tightenumerate}{iii}
\item[\rm(i)] $\Eb[\xi_t \,|\,\Fc_t] = 0$\quad  a.s.;
\item[\rm(ii)]
$\Eb[\|\xi_t\|^2 \,|\,\Fc_t] \le \text{\rm const}$\quad  a.s..
\end{tightenumerate} \vspace{0.5ex}
\end{assumption}

First, we establish an important implication of the ergodicity of the chain. We write $e_i$ for the $i$th unit vector in $\Rb^n$.
\begin{lemma}
\label{l:Poisson}
If the chain $\{i_t\}$ is ergodic with stationary distribution $q$ and Assumption \ref{a:gamma} is satisfied, then
\begin{equation}
\label{error-growth}
 \lim_{T\to \infty}\frac{ \sum_{t=0}^T \gamma_t (e_{i_t} - q)}
 {\sum_{t=0}^T \gamma_t} =0, \quad \text{a.s.},
\end{equation}
and for any $\varepsilon>0$,
\begin{equation}
\label{A3-(b)}
  \lim_{\raisebox{-0.55ex}{${\scriptstyle t_0\to \infty}$}}\sup_{T\ge t_0}\;\frac{ \sum_{t=t_0}^T \gamma_t (e_{i_t} - q)}
 {\max\Big(\varepsilon, \sum_{t=t_0}^T \gamma_t\Big)} =0, \quad \text{a.s.}.
\end{equation}
\end{lemma}
\begin{proof}
Due to the ergodicity of the chain, the vectors
\[
\nu(i) = \Eb \left[ \sum_{t=0}^\infty (e_{i_t}-q)\,\Big| \, i_0=i\right],\quad i\in \Xc,
\]
are finite and satisfy the \emph{Poisson equation}
\begin{equation}
\label{Poisson}
\nu(i) = e_i - q + \sum_{j\in \Xc} P_{ij} \nu(j), \quad i\in \Xc.
\end{equation}
Consider the sums
$\sum_{t=0}^T \gamma_t (e_{i_t} - q)$.
By the Poisson equation,
\begin{equation}
\label{Poisson-1}
e_{i_t} - q = \nu(i_t) - \sum_{j\in \Xc} P_{i_t j} \nu(j)
= \big[ \nu(i_t) - \nu(i_{t+1})\big]  + \Big[ \nu(i_{t+1}) - \sum_{j\in \Xc} P_{i_t j} \nu(j)\Big].
\end{equation}
We consider the two components of the right hand side of \eqref{Poisson-1}, marked with brackets, separately.
Due to Assumption \ref{a:gamma}, (i)---(iii), the series
\[
\sum_{t=1}^\infty \gamma_t \Big[\nu(i_{t+1}) - \sum_{j\in \Xc} P_{i_t j} \nu(j)\Big] =
\sum_{t=1}^\infty \gamma_t \big(\nu(i_{t+1}) - \Eb[\nu(i_{t+1})\,|\,\Fc_t] \big)
\]
is a convergent martingale. Therefore,
\begin{equation}
\notag
 \lim_{T\to \infty}\frac{ \sum_{t=0}^T \gamma_t \big(\nu(i_{t+1}) - \Eb_t[\nu(i_{t+1})] \big)}
 {\sum_{t=0}^T \gamma_t} =0, \quad \text{a.s.}
 \end{equation}
We now focus on the sums
\[
\sum_{t=0}^T \gamma_t \big[\nu(i_t) - \nu(i_{t+1})\big] = \gamma_0 \nu(i_0) + \sum_{t=1}^T (\gamma_t-\gamma_{t-1})\nu(i_t) - \gamma_{T}\nu(i_{T+1}) .
\]
Using Assumption \ref{a:gamma}(iv)  and \cite[Lem. A.3]{ruszczynski1983stochastic},
we obtain \eqref{error-growth}--\eqref{A3-(b)}.
\end{proof}
We can now prove the convergence of the method.
\begin{theorem}
\label{t:TD-convergence}
Suppose the random estimates $\widetilde{\sigma}_{i_t}(P_{i_t},\varPhi r_t)$
 satisfy \eqref{sigma-tilde}, Assumptions \ref{a:gamma} and \ref{a:xi} are satisfied, and  $\alpha\sqrt{1+\varkappa} < 1$.
  If the sequence $\{r_t\}$ is bounded with probability 1,
 then every accumulation point of the sequence $\{r_t\}$ is an element of $Y^*$, with probability 1.
\end{theorem}
\begin{proof}
We use the global Lyapunov function:
\begin{equation}
\label{merit-global}
W(r) =  \min_{r^*\in Y^*}\|   r -  r^*\|^2.
\end{equation}
The direction used in \eqref{TD0-2} at step $t$ can be represented as
\begin{equation}
\label{dir-decom}
\varphi(i_t) \widetilde{d}_t = U(r_t) + \Delta_t,
\end{equation}
with the operator $U(\cdot)$ defined in \eqref{U-operator}, and
\begin{equation}
\label{Delta2}
\Delta_t = -\alpha \xi_t \varphi(i_t) + \varPhi^{\top} \diag\big(e_{i_t}-q\big) \big[ \varPhi r_t - c - \alpha \sigma( P, \varPhi r_t)\big].
\end{equation}
Our intention is to verify the conditions
of Theorem \ref{t:nurminski} for almost all paths of the sequence $\{r_t\}$.
For this purpose, we estimate the decrease of the function \eqref{merit-global}
in iteration $t$. For any $r^*\in Y^*$ we have:
\[
\|r_{t+1}-r^*\|^2 = \|r_t - \gamma_t U(r_t) - r^*\|^2
- 2 \gamma_t \langle \Delta_t, r_t - \gamma_t U(r_t) - r^* \rangle
 + \gamma_t^2 \| \Delta_t\|^2.
\]
The term involving $U(r_t)$ was estimated in the derivation of \eqref{descent-det}. We obtain the
inequality
\begin{equation}
\label{norm-decrease}
\|r_{t+1}-r^*\|^2 \le  \|r_t - r^*\|^2  -  2\gamma_t (1 -\alpha \sqrt{1+\varkappa})  \| \varPhi({r}_{t}-r^*)\|_q^2
- 2 \gamma_t \langle \Delta_t , r_t - \gamma_t U(r_t) - r^* \rangle + C \gamma_t^2 .
\end{equation}

%
Now we can verify the conditions of Theorem \ref{t:nurminski} for almost all paths of the sequence $\{r_t\}$. \\
\emph{Condition A.}  Due to the boundedness of $\{r_t\}$ the sequence $\{U(r_t)\}$ is bounded as well. In view of \eqref{dir-decom},
it is sufficient to verify that $\gamma_t\xi_t\to 0$. By Assumption \ref{a:xi}(i), the sequence
\begin{equation}
\label{xi-martingale}
S_T = \sum_{t=0}^T \gamma_t \xi_t, \quad T=0,1,2,\dots,
\end{equation}
 is a martingale.
Due to  Assumption \ref{a:xi}(ii), $\Eb[S_T^2] \le \text{const}\cdot \Eb[\sum_{t=0}^T\gamma_t^2]$. In view of  Assumption \ref{a:gamma}(ii), by virtue of the martingale convergence theorem, $\{S_T\}$ is convergent a.s.,
which yields $\lim_{t\to \infty} \gamma_t \xi_t = 0$. \\
\emph{Condition B.}
Suppose  $r_k\to r' \notin Y^*$ for $k\in \Kc$ (on a certain path $\omega$).
If B were false, then for all $\varepsilon_0>0$ we could find $\varepsilon\in (0,\varepsilon_0]$ and $k\in \Kc$ such that
$\|r_t-r_k\| \le \varepsilon$ for all $t \ge k$. Then for all $k_0\in \Kc$, $k_0 \ge k$, we have $\|r_t-r_{k_0}\| \le 2\varepsilon$ for all $t \ge k_0$. Since $r'$ is not optimal,
we can choose $\varepsilon_0>0$ small enough, $k_0\in \Kc$ large enough, and $\delta>0$ small enough, so that $\|\varPhi(r_t-r^*)\|^2_q > \delta$ for all $t \ge k_0$. Then
\eqref{norm-decrease} yields
\begin{multline}
\label{multi-decrease}
\|r_{T}-r^*\|^2 \le  \|r_{k_0} - r^*\|^2\\
{} + \left(  -\delta(1-\alpha\sqrt{1+\varkappa})
+ \frac{\sum_{t=k_0}^{T-1}\gamma_t \langle \Delta_t, r_t - \gamma_t U(r_t) - r^*\rangle}{\sum_{t=k_0}^{T-1} \gamma_t } + C\frac{\sum_{t=k_0}^{T-1}\gamma_t^2}{\sum_{t=k_0}^{T-1} \gamma_t }\right)\sum_{t=k_0}^{T-1} \gamma_t.
\end{multline}
We fix $r^*=\text{Proj}_{Y^*}(r_{k_0})$ and estimate the growth of the sums involving $\Delta_t$. We write $\Delta_t=\Delta^{(1)}_t+\Delta^{(2)}_t$,
where, in view of \eqref{Delta2},
\[
\Delta^{(1)}_t = -\alpha \xi_t \varphi(i_t), \quad \Delta^{(2)}_t=\varPhi^{\top} \diag\big(e_{i_t}-q\big) \big[ \varPhi r_t - c - \alpha \sigma( P, \varPhi r_t)\big].
\]

Since \eqref{xi-martingale} is a convergent martingale and the terms $ \langle \varphi(i_t), r_t - \gamma_t U(r_t) - r^*\rangle$ are bounded and $\Fc_t$-measurable, we have
\[
\lim_{T\to\infty}\frac{\left|\sum_{t=k_0}^{T-1}\gamma_t  \langle \Delta^{(1)}_t, r_t - \gamma_t U(r_t) - r^*\rangle\right|}{\sum_{t=k_0}^{T-1} \gamma_t } = 0.
\]
To deal with the sum involving $\Delta^{(2)}_t$, observe that $\| r_t-r_{k_0}\|\le 2 \varepsilon_0$ and thus
\begin{multline}
\label{Delta2b}
\langle \Delta^{(2)}_t,r_t - \gamma_t U(r_t) - r^*\rangle =\\
 \Big\langle  \diag\big(e_{i_t}-q\big)
\big[ \varPhi r_{k_0} - c - \alpha \sigma( P, \varPhi r_{k_0})\big],\varPhi (r_{k_0}-r^*)\Big\rangle + h_t
= \langle e_{i_t}-q , w\rangle + h_t,
\end{multline}
where $|h_t| \le C\varepsilon_0$ and $w$ is a fixed vector (depending on $k_0$ only). It follows that
\begin{equation}
\label{Delta2-series}
\left| \sum_{t=k_0}^{T-1}\gamma_t \langle \Delta^{(2)}_t, r_t - \gamma_t U(r_t) - r^*\rangle \right| \le
C \left\|\sum_{t=k_0}^{T-1} \gamma_t (e_{i_t} - q)\right\| +C \varepsilon_0 \sum_{t=k_0}^{T-1} \gamma_t.
\end{equation}
Dividing both sides of \eqref{Delta2-series} by $\sum_{t=k_0}^{T-1} \gamma_t$ and using \eqref{error-growth}, we see that we can choose $\varepsilon>0$ small enough and $k_0\in \Kc$ large enough, so that the entire expression
in parentheses in \eqref{multi-decrease} is smaller than $-\delta(1-\alpha\sqrt{1+\varkappa})/2$, if $T$ is large enough. But this yields $\|r_T-r^*\| \to -\infty$, as $T\to\infty$, a contradiction.
Therefore, Condition B is satisfied.\\
\emph{Condition C.} The inequality \eqref{multi-decrease} remains valid for $T=s(k_0,\varepsilon)$. By the definition
of $s(k_0,\varepsilon)$,
\[
\Big\| \sum_{t=k_0}^{T-1} \gamma_t(d_t+\xi_t)\Big\| \ge \varepsilon.
\]
By the convergence of
\eqref{xi-martingale}, and the boundedness of $\{d_t\}$, a constant $C>0$ exists  such that for all sufficiently large $k_0$ and sufficiently small $\varepsilon$, we have
\[
\sum_{t=k_0}^{T-1} \gamma_t \ge \varepsilon/C.
\]
Using \eqref{A3-(b)}, by a similar argument as in the analysis of Condition B,  we can
choose $\varepsilon_1>0$ small enough that for all $k_0\in \Kc$ large enough so that the entire expression
in parentheses in \eqref{multi-decrease} is smaller than $-\delta(1-\alpha\sqrt{1+\varkappa})/2$. Therefore, for all $\varepsilon\in (0,\varepsilon_1]$ and all sufficiently large $k_0\in \Kc$
\[
\|r_{s(k_0,\varepsilon)}-r^*\|^2  \le \|r_{k_0}- r^*\|^2 - \frac{\delta(1-\alpha\sqrt{1+\varkappa})\varepsilon}{2C}.
\]
We fix $r^*=\text{Proj}_{Y^*}(r_{k_0})$ on the right hand side, and  obtain
\[
W(r_{s(k_0,\varepsilon)}) \le  \|r_{s(k_0,\varepsilon)}-r^*\|^2\le   W(r_{k_0}) - \frac{\delta(1-\alpha\sqrt{1+\varkappa})\varepsilon}{2C}.
\]
Now, the limit with respect to $k_0\to \infty$, $k_0\in \Kc$, proves
Condition C.\\
\emph{Condition D} is satisfied trivially, because $W(r^*)\equiv 0$ for $r^*\in Y^*$.
\end{proof}

The only question remaining is the boundedness of the sequence $\{r_t\}$. It is a common issue in the analysis of stochastic approximation
algorithms \cite[\S 5.1]{kushner2003stochastic}. In our case,  no additional conditions and analysis are needed, because our Lyapunov function
\eqref{merit-global} is the squared distance to the optimal set. Therefore, a simple algorithmic modification: the projection on a bounded set $Y$ intersecting with $\{r \in \Rb^m: \varPhi r = v^*\}$, is sufficient to guarantee boundedness. The modified method \eqref{TD0-2} reads:
\begin{equation}
\label{TD0-proj}
r_{t+1} =  \proj_{Y}\big( r_t - \gamma_t \varphi(i_t)\, \widetilde{d}_t\big), \quad t=0,1,2,\dots.
\end{equation}
Now,
$Y^* = \{r \in Y: \varPhi r = v^*\}$  and we require that this set is nonempty.
This modification does not affect our analysis in any meaningful way, because the projection is nonexpansive.
In the proof of Theorem \ref{t:RTD-convergence}, we use the inequality
\[
\big\|\proj_{Y} (r' - \gamma U(r'))- \proj_{Y}(r'' - \gamma U(r''))\big\|^2 \le
\big\| (r' - \gamma U(r'))- (r'' - \gamma U(r''))\big\|^2
\]
and proceed as before. In the proof of Theorem \ref{t:TD-convergence}, we start from
\[
\|r_{t+1}-r^*\|^2 = \big\|\proj_Y\big(r_t - \gamma_t (U(r_t)+\Delta_t)\big) - r^*\big \|^2 \le \big\|r_t - \gamma_t (U(r_t)+\Delta_t) - r^*\big \|^2,
\]
and then continue in the same way as before. We did not include projection into the method originally, because it obscures
the presentation. In practice, we have not yet encountered any need for it.

\section{The Multistep Risk-Averse Method of Temporal Differences}
\label{s:TDL}

In the method  discussed so far, the residuals are corrected by moving in the direction of the last feature vector $\varphi(i_t)$. Alternatively,  we may use the weighted averages of all previous observations, where the highest weight is given to the most recent observation and the weights decrease exponentially as we look into the past observations. This idea is the core of the well-known TD($\lambda$) algorithm \cite{Sutton1988}. We generalize it to the risk-averse case.

For a fixed policy $\pi$,  we refer to $v^{\pi}$ as $v$, and to $P^{\pi}$ as $P$, for simplicity. The multistep risk-averse method of temporal differences carries out the following iterations:
\begin{gather}
z_t = \lambda \alpha z_{t-1} + \varphi(i_t),\quad t = 0,1,2,\ldots, \label{z-update}\\
r_{t+1} = r_{t} - \gamma_t z_t \widetilde{d}_t,\quad t =0, 1,2,\ldots \label{r-update}
\end{gather}
where $\lambda\in [0,1]$, and $\widetilde{d}_t$ is given by \eqref{TD0-1}.
For simplicity, $z_{-1}$ is assumed to be the zero vector. In the risk-neutral case, when
$\sigma_{i_t} (P_{i_t},\varPhi r_t) = P_{i_t} \varPhi r_t$, the method reduces to the
classical TD($\lambda$).

Our convergence analysis will use some ideas from the analysis in the previous two sections, albeit in
a form adapted to the version with exponentially averaged features. However, contrary to the
expected value setting,  the
method \eqref{z-update}--\eqref{r-update}  will converge to a solution of an equation different from \eqref{DP-projected},
but still relevant for our problem.

We start from a heuristic analysis of a deterministic
counterpart of the method, to extract its  drift. In the next section, we make all approximations precise, but we believe that this introduction is useful to decipher our detailed approach to follow.
By direct calculation,
\begin{equation}
\label{zt-expansion}
 z_t =  \sum_{k=0}^t (\lambda\alpha)^{t-k} \varphi(i_k),
\end{equation}
and thus
\[
 z_t d_t = \varPhi^\top \sum_{k=0}^t (\lambda\alpha)^{t-k}  e_{i_k} e_{i_t}^\top  \big( \varPhi r_t - c - \alpha {\sigma} (P,\varPhi r_t)\big).
\]
Heuristically assuming that $r_t \approx r'$, we focus on the operator acting on the expected temporal differences. As each of the observed
feature vectors $\varphi(i_k)$ affects all succeeding steps of the method, via the filter \eqref{z-update}, we need to study the cumulative effect of many steps. We look, therefore, at the sums
\[
G_T = \Eb\bigg[ \sum_{t=0}^T \gamma_t \sum_{k=0}^t (\lambda\alpha)^{t-k}  e_{i_k} e_{i_t}^\top \bigg].
\]
Changing the order of summation and using the fact that $\{(\lambda\alpha)^{t-k}\}_{t \ge k}$ diminishes very fast, as compared to
$\{\gamma_t\}_{t \ge k}$,   we get
\[
G_T =  \Eb\bigg[ \sum_{k=0}^T \sum_{t=k}^T \gamma_t  (\lambda\alpha)^{t-k}  e_{i_k} e_{i_t}^\top \bigg] \approx
\Eb\bigg[ \sum_{k=0}^T \gamma_k \sum_{t=k}^T (\lambda\alpha)^{t-k}  e_{i_k} e_{i_t}^\top \bigg].
\]
Therefore
\begin{align*}
G_T &\approx
\Eb\bigg[ \sum_{k=0}^T \gamma_k \sum_{t=k}^T   (\lambda\alpha)^{t-k}  e_{i_k}\Eb\big[  e_{i_t}^\top \,\big|\,\Fc_k\big] \bigg]
=
\Eb\bigg[ \sum_{k=0}^T \gamma_k \sum_{t=k}^T   (\lambda\alpha)^{t-k}  e_{i_k}e_{i_k}^\top P^{t-k} \bigg]\\
&=  \sum_{k=0}^T \gamma_k \Eb[ \diag(e_{i_k})\big] \sum_{t=k}^T   (\lambda\alpha)^{t-k}   P^{t-k}
\approx \sum_{k=0}^T \gamma_k \Eb[ \diag(e_{i_k})\big] \sum_{t=k}^\infty  (\lambda\alpha)^{t-k}   P^{t-k}\\
&\approx Q \sum_{k=0}^T \gamma_k \sum_{t=k}^\infty  (\lambda\alpha)^{t-k}   P^{t-k} .
\end{align*}
The last approximations are possible because $\lambda\alpha\in [0,1)$ and $\Eb[ \diag(e_{i_k})\big] \to q$ at an exponential rate.
We now define the multistep transition matrix,
\begin{equation}
\label{bar-P}
\widebar{P}  = (1-\lambda\alpha)\sum_{\ell=0}^\infty  (\lambda\alpha)^{\ell}P^{\ell}.
\end{equation}
By construction, $\widebar{P}\in \widebar{\conv}\{I,P,P^2,\dots\}$.
With these approximations, we can simply write
\[
G_T \approx \frac{1}{1-\lambda\alpha} Q  \widebar{P} \sum_{k=0}^T {\gamma}_k.
\]
Define the operators
\begin{equation}
\label{bar-U-operator}
\widebar{U}(r)  =   \varPhi^{\top} {Q} \widebar{P} \big[ \varPhi r -c - \alpha \,\sigma(P,\varPhi r  )\big], \quad t=0,1,2,\dots
\end{equation}
and consider the following deterministic counterpart of \eqref{z-update}--\eqref{r-update}, with $\widebar{\gamma} \sim {\gamma_t}/{(1-\lambda\alpha)}$:
\begin{equation}
\label{RTDL-det}
r_{r+1} = r_t - \widebar{\gamma}\, \widebar{U}(r_t),\quad t=0,1,2,\dots, \quad \widebar{\gamma}>0.
\end{equation}
Our intention is to show that for sufficiently small $\widebar{\gamma}$ the method \eqref{RTDL-det} converges to a point $r^*$ such that
$\widebar{U}(r^*)=0$. Such a point is also a solution of the following
projected multistep risk-averse dynamic programming equation:
\begin{align}
  L\widebar{P}\varPhi r  =   L\widebar{P}\big(  c + \alpha \sigma(P,\varPhi r)\big), \label{projeq-L}
\end{align}
where $L$ is the projection operator defined in \eqref{ldef}.
The solutions of \eqref{projeq-L}  differ from the solutions
of \eqref{projeq}, unlike in the risk-neutral case
($\varkappa=0$). If we replace $\widebar{P}$ with  $I$, \eqref{projeq-L} reduces to \eqref{projeq}.


\begin{theorem}
\label{RTDL-convergence}
If $\alpha (1+\varkappa)< 1$, then $\widebar{\gamma}_0>0$ exists, such that for all $\widebar{\gamma}\in (0,\widebar{\gamma}_0)$
the algorithm \eqref{RTDL-det} generates a sequence $\{r_t\}$ convergent to a point $r^*$ such that $\widebar{U}(r^*)=0$.
\end{theorem}
\begin{proof}
For two arbitrary points $r'$ and $r''$ we have
\begin{multline}
\label{step-expansion}
\Big\| \big(r' - \widebar{\gamma}\, \widebar{U}(r')\big)- \big(r'' - \widebar{\gamma}\, \widebar{U}(r'')\big)\Big\|^2
\\=  \|r'-r''\|^2 \
\ + 2 \widebar{\gamma}\, \Big\langle r' - r'' ,  \varPhi^{\top}Q\widebar{P} \big[ - \varPhi (r' - r'')
  +\alpha  \sigma(P,\varPhi r')- \alpha \sigma(P,\varPhi r'')\big]\Big\rangle\\
{ } + \widebar{\gamma}^2 \Big\| \varPhi^{\top} Q \widebar{P} \varPhi (r' - r'')  - \alpha\varPhi^{\top} Q \widebar{P}\big[ \sigma(P,\varPhi r')-\sigma(P,\varPhi r'')\big]\Big\|^2.
\end{multline}
We focus on the scalar product in the middle of the right hand side of \eqref{step-expansion}:
\begin{multline}
\label{scalar-Roy}
\Big\langle \varPhi(r' - r'') ,  \widebar{P} \big[ - \varPhi (r' - r'')
  +\alpha  \sigma(P,\varPhi r')- \alpha \sigma(P,\varPhi r'')\big]\Big\rangle_q\\
  = \Big\langle \varPhi(r' - r'') ,  \widebar{P} \big[ - \varPhi (r' - r'') +\alpha P \varPhi (r' - r'') \big]\Big\rangle_q \\
  {} +
 \alpha \Big\langle \varPhi(r' - r'') ,  \widebar{P} \big[
    \sigma(P,\varPhi r')-  \sigma(P,\varPhi r'')-  P \varPhi (r' - r'')  \big]\Big\rangle_q .
\end{multline}
Setting $h=\varPhi (r' - r'')$, we can estimate the first (quadratic) term on the right hand side of \eqref{scalar-Roy} by a calculation
borrowed from \cite[Lem. 8]{tsitsiklis1997analysis}, with $h=\varPhi (r' - r'')$:
\begin{align*}
\lefteqn{\Big\langle  h,  \widebar{P} \big[ - h +\alpha Ph \big]\Big\rangle_q =(1-\alpha\lambda) \Big\langle  h, \sum_{\ell=0}^\infty  (\alpha\lambda)^{\ell}P^{\ell} \big[ - h +\alpha P h \big]\Big\rangle_q}\qquad\\
&=(1-\alpha\lambda)(1-\lambda) \Big\langle  h, \sum_{k=0}^\infty \lambda^k \sum_{\ell=0}^k  \alpha^{\ell}P^{\ell} \big[ - h +\alpha P h \big]\Big\rangle_q\\
&=(1-\alpha\lambda)(1-\lambda) \Big\langle  h, \sum_{k=0}^\infty \lambda^k  \big[ \alpha^{k+1}P^{k+1}h  - h \big]\Big\rangle_q\\
&=(1-\alpha\lambda) \Big\langle  h, (1-\lambda)\sum_{k=0}^\infty \lambda^k   \alpha^{k+1}P^{k+1}h  - h \Big\rangle_q\\
&=(1-\alpha\lambda) \Big\langle  h, \frac{\alpha(1-\lambda)}{1- \alpha\lambda} P \widebar{P}h  - h \Big\rangle_q
\le  (\alpha - 1) \|h\|_q^2.
\end{align*}
The last inequality is due to the fact that both $P$ and $\widebar{P}$ are nonexpansive in $\|\cdot\|_q$.

The second (nonsmooth) term on the right hand side of \eqref{scalar-Roy} can be estimated by \eqref{sigma-delta}, again with the use of the nonexpansiveness of  $\widebar{P}$:
\[
\Big\langle \varPhi(r' - r'') ,  \widebar{P} \big[
    \sigma(P,\varPhi r')-  \sigma(P,\varPhi r'')-  P \varPhi (r' - r'')  \big]\Big\rangle_q \le \varkappa \,\big\| \varPhi (r' - r'')\big\|_q^2.
\]
The last term on the right hand side of \eqref{step-expansion} (with ${\gamma}^2$) can be bounded by
 $\gamma^2 \widebar{C} \|\varPhi (r'-r'')\|_q^2$, where $\widebar{C}$ is some constant.
Integrating all these estimates into \eqref{step-expansion}, we obtain the inequality
\[
\big\| (I - \widebar{\gamma}\, \widebar{U})(r')- (I - \widebar{\gamma}\, \widebar{U})(r'')\big\|^2
 \le  \|r'-r''\|^2 - 2 \widebar{\gamma}\, \Big(1 - \alpha(1 + \varkappa) - \frac{\widebar{\gamma}\, \widebar{C}}{2}\Big)\big\| \varPhi (r' - r'')\big\|_q^2 .
\]
If $\alpha (1+\varkappa)< 1$, then using $0< \widebar{\gamma} < 2(1- \alpha(1+\varkappa))/\widebar{C} $, we obtain:
\begin{equation}
\label{semi-contraction-L}
\big\| (I - \widebar{\gamma}\, \widebar{U})(r')- (I - \widebar{\gamma}\, \widebar{U})(r'')\big\|^2  \le \|r'-r''\|^2 - \widebar{\gamma} \beta \|\varPhi (r'-r'')\|_q^2,
\end{equation}
with some $\beta>0$. In particular, setting $r' = \bar{r}_t$ and $r''=r^*$ for a solution $r^*$ of \eqref{projeq}, we obtain
the following relation between successive iterates of the method \eqref{RTDL-det}:
\begin{equation}
\label{descent-det-L}
\| {r}_{t+1}-r^*\|^2 \le \| {r}_{t}-r^*\|^2 - \widebar{\gamma} \beta  \| \varPhi({r}_{t}-r^*)\|_q^2.
\end{equation}
This immediately proves that the sequence $\{{r}_t\}$ is bounded and $\varPhi {r}_t \to \varPhi r^*$. Every accumulation point $\hat{r}$ of $\{{r}_t\}$ must be then a solution
of equation \eqref{projeq-L}. Substituting this accumulation point for $r^*$ in the last inequality, we conclude that the entire sequence $\{{r}_t\}$
is convergent to $\hat{r}$.
\end{proof}

If $\varPhi$ has full column rank, the solution $r^*$ is unique, because substituting another solution for $r_t$ in \eqref{RTDL-det} we obtain $r_{t+1} = r_t$, which
leads to a contradiction in \eqref{descent-det-L}.

\section{Convergence of the Risk-Averse Multistep Method}
\label{s:TDL-convergence}

We now carry out a detailed analysis of the stochastic method \eqref{z-update}--\eqref{r-update}.

\begin{lemma}
\label{l:gamma-simplification}
For any array of uniformly bounded random variables $\big\{ A_{k,t}\big\}_{k \ge 0,\, t\ge 0} $
\[
\lim_{T\to \infty} \frac{ \sum_{k=0}^T\sum_{t=k}^T \gamma_t (\lambda\alpha)^{t-k} A_{k,t} - \sum_{k=0}^T\gamma_k \sum_{t=k}^\infty  (\lambda\alpha)^{t-k} A_{k,t}}
{ \sum_{k=0}^T \gamma_k} = 0, \quad  \text{a.s.}
\]
\end{lemma}
\begin{proof}
Changing the order of summation twice, we obtain
\begin{align*}
\lefteqn{ \sum_{k=0}^T\sum_{t=k+1}^T |\gamma_t - \gamma_k| (\lambda\alpha)^{t-k} \le
 \sum_{k=0}^T\sum_{t=k+1}^T \sum_{\ell=k+1}^t|\gamma_{\ell} - \gamma_{\ell-1}| (\lambda\alpha)^{t-k}}\\
& \le \frac{1}{1-\lambda\alpha}\sum_{k=0}^T\sum_{\ell=k+1}^T |\gamma_{\ell} - \gamma_{\ell-1}|(\lambda\alpha)^{\ell -k}
 = \frac{1}{1-\lambda\alpha} \sum_{\ell=1}^T|\gamma_{\ell} - \gamma_{\ell-1}|\sum_{k=0}^{\ell-1} (\lambda\alpha)^{\ell -k}\\
& \le \frac{\lambda\alpha}{(1-\lambda\alpha)^2} \sum_{\ell=1}^T|\gamma_{\ell} - \gamma_{\ell-1}|.
\end{align*}
Therefore, with $C$ being the uniform bound on $\|A_{k,t}\|$ and $\gamma_{k}^{\,\max} = \max_{t \ge k} \gamma_t$, we obtain
\begin{align*}
\lefteqn{\bigg\| \sum_{k=0}^T\sum_{t=k}^T \gamma_t (\lambda\alpha)^{t-k} A_{k,t} - \sum_{k=0}^T\gamma_k \sum_{t=k}^\infty  (\lambda\alpha)^{t-k} A_{k,t} \bigg\|}\\
&\le C \sum_{k=0}^T\sum_{t=k+1}^T |\gamma_t - \gamma_k| (\lambda\alpha)^{t-k} + C \sum_{k=0}^T\sum_{t=T+1}^\infty \gamma_t (\lambda\alpha)^{t-k}\\
 &\le \frac{C\lambda\alpha}{(1-\lambda\alpha)^2} \sum_{\ell=1}^T|\gamma_{\ell} - \gamma_{\ell-1}| +  \frac{C \gamma_{T+1}^{\,\max}\lambda\alpha}{(1-\lambda\alpha)^2}.
\end{align*}
Assumption \ref{a:gamma}(iv) and \cite[Lem. A.3]{ruszczynski1983stochastic} imply the assertion.
\end{proof}

We need another auxiliary result, extending Lemma \ref{l:Poisson} to our case.
\begin{lemma}
\label{l:Poisson-lambda}
 \begin{equation}
 \label{ST-growth}
 \lim_{T\to \infty} \frac
 {{\sum_{t=0}^T \gamma_t \bigg( \sum_{k=0}^t  (\lambda\alpha)^{t-k} e_{i_k}e_{i_t}^\top- \frac{1}{1-\lambda\alpha} Q \widebar{P}\bigg)}}
 {{\sum_{t=0}^T \gamma_t}}  = 0 \quad \text{a.s.},
 \end{equation}
 and for any $\varepsilon>0$,
\begin{equation}
\label{ST-growth-b}
  \lim_{\raisebox{-0.55ex}{${\scriptstyle t_0\to \infty}$}}\sup_{T\ge t_0}\;
  \frac{ \sum_{t=t_0}^T \gamma_t \bigg( \sum_{k=0}^t  (\lambda\alpha)^{t-k} e_{i_k}e_{i_t}^\top- \frac{1}{1-\lambda\alpha} Q \widebar{P}\bigg)}
 {\max\Big(\varepsilon, \sum_{t=t_0}^T \gamma_t\Big)} =0  \quad \text{a.s.}.
\end{equation}
 \end{lemma}
 \begin{proof} Consider the sums appearing in the numerator of \eqref{l:Poisson-lambda}:
 \[
\sum_{t=0}^T \gamma_t \sum_{k=0}^t  (\lambda\alpha)^{t-k} e_{i_k}e_{i_t}^\top
= \sum_{k=0}^T e_{i_k}\sum_{t=k}^T (\lambda\alpha)^{t-k}\gamma_t e_{i_t}^\top,\quad T=1,2,\dots.
 \]
In view of Lemma \ref{l:gamma-simplification}, it is sufficient to consider the sums
\[
S_T = \sum_{k=0}^T \gamma_k\, e_{i_k}\sum_{t=k}^\infty (\lambda\alpha)^{t-k} e_{i_t}^\top,\quad T=1,2,\dots.
\]
We transform the inner sum:
 \begin{multline*}
\sum_{t=k}^\infty ((\lambda\alpha))^{t-k} e_{i_t}^\top = \sum_{t=k}^\infty (\lambda\alpha)^{t-k} \Big\{
\sum_{\ell=k+1}^t  \big[  e_{i_\ell}^\top P^{t-\ell}  - e_{i_{\ell-1}}^\top P^{t-\ell+1} \big] +   e_{i_{k}}^\top P^{t-k} \Big\} \\
= \sum_{t=k}^\infty  (\lambda\alpha)^{t-k} e_{i_{k}}^\top P^{t-k}
+ \sum_{\ell=k+1}^\infty \sum_{t=\ell}^\infty  (\lambda\alpha)^{t-k}
 \big[  e_{i_\ell}^\top P^{t-\ell}  - e_{i_{\ell-1}}^\top P^{t-\ell+1} \big].
\end{multline*}
We can thus write $S_T = S^{(1)}_T+S^{(2)}_T$, with
\[
S^{(1)}_T =  \sum_{k=0}^T \gamma_k\, e_{i_k}e_{i_{k}}^\top \sum_{t=k}^\infty  (\lambda\alpha)^{t-k}  P^{t-k}  =
\frac{1}{1-\lambda\alpha}\sum_{k=0}^T \gamma_k\, \diag(e_{i_k}) \widebar{P}
\]
and
\begin{align*}
S^{(2)}_T  &= \sum_{k=0}^T\gamma_k\, e_{i_k} \sum_{\ell=k+1}^\infty \sum_{t=\ell}^\infty   (\lambda\alpha)^{t-k}
 \big[  e_{i_\ell}^\top P^{t-\ell}  - e_{i_{\ell-1}}^\top P^{t-\ell+1} \big] \\
 &= \frac{1}{1-\lambda\alpha}\sum_{k=0}^T\gamma_k\, e_{i_k} \sum_{\ell=k+1}^\infty (\lambda\alpha)^{\ell-k}
 \big[  e_{i_\ell}^\top  - e_{i_{\ell-1}}^\top P \big] \widebar{P}\\
 &= \frac{1}{1-\lambda\alpha}\sum_{\ell=1}^\infty\sum_{k=0}^{\min(T,\ell-1)}\gamma_k\, e_{i_k} (\lambda\alpha)^{\ell-k} \big[  e_{i_\ell}^\top  - e_{i_{\ell-1}}^\top P \big] \widebar{P}.
\end{align*}
The second sum is a convergent martingale, because
$\Eb\big[e_{i\ell}^\top\,\big|\,\Fc_{\ell-1}\big] = e_{i_{\ell-1}}^\top P$.
Therefore, it satisfies \eqref{ST-growth}.

Applying Lemma \ref{l:Poisson} to  $S^{(1)}_T  - \frac{1}{1-\lambda\alpha}\sum_{k=0}^T \gamma_k\, \diag(q) \widebar{P}$,
we obtain both assertions.
 \end{proof}

Now we can follow the arguments of \S \ref{s:TD0-convergence} and establish the convergence of the
multistep method.

\begin{theorem}
\label{t:TDL-convergence}
Assume that $\alpha(1+\varkappa) < 1$, the sequence $\{r_t\}$ is bounded with probability~1,
 and
the random estimates $\widetilde{\sigma}_{i_t}(P_{i_t},\varPhi r_t)$
 satisfy \eqref{sigma-tilde}.
 Then, with probability~1, every accumulation point of the sequence $\{r_t\}$ generated by
 \eqref{z-update}--\eqref{r-update}  is a solution of \eqref{projeq-L}.
\end{theorem}
\begin{proof}
We represent the direction used in \eqref{r-update} at step $t$ as
\[
z_t \widetilde{d}_t =  \frac{1}{1-\lambda\alpha}\widebar{U}_t(r_t) + \Delta^{(1)}_t + \Delta^{(2)}_t,
\]
with the operator $\widebar{U}_t(\cdot)$ defined in \eqref{bar-U-operator},
and
\begin{align*}
\Delta^{(1)}_t &= -\alpha z_t \xi_t,\\
\Delta^{(2)}_t & = z_t d_t -  \frac{1}{1-\lambda\alpha}\widebar{U}_t(r_t).
\end{align*}
For any $r^*$ solving \eqref{projeq-L}, with $\widebar{\gamma}_t = \gamma_t/(1-\lambda\alpha)$, we have
\[
\big\|r_{t+1}-r^* \big\|^2 = \big\|r_t - \widebar{\gamma}_t  \widebar{U}_t(r_t) - r^*\big\|^2
- 2 \widebar{\gamma}_t \langle \Delta^{(1)}_t + \Delta^{(2)}_t, r_t - \widebar{\gamma}_t  \widebar{U}_t(r_t) - r^* \rangle
 + \widebar{\gamma}_t^2 \big\| \Delta^{(1)}_t+\Delta^{(2)}_t \big\|^2.
\]
 Our intention is to verify the conditions
of Theorem \ref{t:nurminski} for almost all paths of the sequence $\{r_t\}$.

\emph{Condition A.}  The sequence $\{z_t\}$ is bounded by construction. Since the series \eqref{xi-martingale} is a convergent martingale,
we conclude that $\lim_{t\to \infty} \gamma_t z_t \widetilde{d}_t = 0$. \\

\emph{Conditions B and C}: We follow the proof of Theorem \ref{t:TD-convergence}.
The deterministic term involving $\widebar{U}_t(r_t)$ can be estimated as in \eqref{descent-det-L}:
\[
\big\|r_t - \widebar{\gamma}_t \widebar{U}_t (r_t) - r^*\big\|^2 \le
\big\|r_t - r^*\big\|^2  -  2\widebar{\gamma}_t  \big(1 -\alpha(1+\varkappa)\big)  \big\| \varPhi({r}_{t}-r^*)\big\|_q^2 + C\widebar{\gamma}_t^2.
\]
Since $\{z_t\}$ and $\{r_t\}$ are bounded,
Assumptions \ref{a:gamma} and \ref{a:xi} imply that
 $
\sum_{t=0}^\infty \widebar{\gamma}_t \langle \Delta^{(1)}_t, r_t - \widebar{\gamma}_t \widebar{U}_t(r_t) - r^*\rangle$
 is a convergent martingale.

To analyze the second error term,  $\Delta^{(2)}_t$,  we observe that for a vector $e_{i_k}$ having 1 at position~$i_k$ and zero otherwise, the formula \eqref{zt-expansion} yields
\begin{align*}
z_td_t &=\sum_{k=0}^t (\lambda\alpha)^{t-k} \varphi(i_k) \big( \varphi^{\top}\!(i_t)r_t - c(i_t) - \alpha {\sigma}_{i_t} (P_{i_t},\varPhi r_t) \big)\\
&=
\varPhi^{\top}\Big( \sum_{k=0}^t (\lambda\alpha)^{t-k}e_{i_k}e_{i_t}^\top\Big) \big( \varPhi r_t -
c - \alpha \sigma( P, \varPhi r_t)\big).
\end{align*}
Subtracting \eqref{bar-U-operator},  we obtain
\[
\Delta^{(2)}_t = \varPhi^{\top} \Big( \sum_{k=0}^t (\lambda\alpha)^{t-k} e_{i_k}e_{i_t}^\top - \frac{1}{1-\lambda\alpha}Q \widebar{P}_t  \Big) \big[ \varPhi r_t - c - \alpha \sigma( P, \varPhi r_t)\big].
\]
By virtue of Lemma \ref{l:Poisson-lambda}, for any $\varepsilon>0$,
\[
\lim_{T\to \infty} \frac{\sum_{t=0}^T \gamma_t \Delta^{(2)}_t}{\sum_{t=0}^T \gamma_t} = 0, \qquad
  \lim_{\raisebox{-0.55ex}{${\scriptstyle t_0\to \infty}$}}\sup_{T\ge t_0}\;
  \frac{ \sum_{t=t_0}^T \gamma_t \Delta^{(2)}_t}
 {\max\Big(\varepsilon, \sum_{t=t_0}^T \gamma_t\Big)} =0  \quad \text{a.s.}.
\]

The remaining analysis is  the same as in the proof of Theorem \ref{t:TD-convergence}.
We obtain an inequality corresponding to \eqref{multi-decrease}:
\begin{multline}
\notag
\|r_{T}-r^*\|^2 \le  \|r_{k_0} - r^*\|^2\\
{} + \left(  -\delta(1-\alpha(1+\varkappa))
+ \frac{\sum_{t=k_0}^{T-1}{\gamma}_t \langle \Delta^{(1)}_t+\Delta^{(2)}_t, r_t - \widebar{\gamma}_t  \widebar{U}(r_t) - r^*\rangle}{\sum_{t=k_0}^{T-1} \widebar{\gamma}_t } + C\frac{\sum_{t=k_0}^{T-1}\widebar{\gamma}_t^2}{\sum_{t=k_0}^{T-1} \widebar{\gamma}_t }\right)\sum_{t=k_0}^{T-1} \widebar{\gamma}_t,
\end{multline}
with $\delta>0$.
This allows us to verify the conditions of Theorem \ref{t:nurminski} and prove our assertion following the last steps of the proof of Theorem  \ref{t:TD-convergence}
{verbatim}.
\end{proof}
It is worth mentioning that the convergence condition  for the multistep method: $\alpha(1+\varkappa) <1$, is slightly stronger that the condition
for the basic method: $\alpha \sqrt{1+\varkappa}<1$.

Again, as in the case of the basic method, discussed in \S \ref{s:TD0-convergence}, the boundedness of the
sequence $\{r_k\}$ is not an issue of concern, because it can be guaranteed by projection on a bounded set $Y$. The modified method has the following form:
\begin{equation}
\label{TDL-proj}
r_{t+1} =  \proj_{Y}\big( r_t - \gamma_t z_t\, \widetilde{d}_t\big), \quad t=0,1,2,\dots.
\end{equation}
We just need $Y$ to have a nonempty intersection $Y^*$ with the set of solutions of \eqref{projeq-L}.
Due to the nonexpansiveness of the projection operator, all our proofs remain unchanged with this modification, as
discussed at the end of \S \ref{s:TD0-convergence}.

\section{Empirical Study}
\label{s:empirical}

\subsection{Risk estimation}

We first discuss the issue of obtaining stochastic estimates  $\widetilde{\sigma}_{i_t}(P_{i_t},\cdot)$ satisfying
\eqref{sigma-tilde} and Assumption \ref{a:xi}:
\begin{equation}
\label{sigma-tilde-E}
\Eb\big[ \widetilde{\sigma}_{i_t}(P_{i_t},\varPhi r_t)\big|\, \Fc_t\big]  = \sigma_{i_t} (P_{i_t},\varPhi r_t),\quad t=0,1,2,\dots,
\end{equation}
In the expected value case, where
$\sigma_{i_t} (P_{i_t},\varPhi r_t) = P_{i_t} \varPhi r_t = \Eb\big[\varphi^\top(i_{t+1}) r_t\,\big|\,\Fc_t\big]$,
we could just use the approximation value  at the next state observed, $\varphi^\top(i_{t+1}) r_t$, as the stochastic estimate of the
expected value function. However, due to the nonlinearity of a risk measure with respect to the probability measure
$P_{i_t}$, such a straightforward approach is no longer possible.

Statistical estimation of measures of risk is a challenging problem, for which, so far, only solutions in special cases
have been found \cite{dentcheva2017statistical}. To mitigate this problem, we propose to use a special class of transition risk mappings which are very convenient
for statistical estimation. For a given transition risk mapping $\widebar{\sigma}_i(P_i,v)$, we sample $N$ conditionally independent transitions
from the state $i$, resulting in states $j^1, \dots,j^N$. This sample defines a random empirical distribution,
$
P_i^N  = \frac{1}{N} \sum_{k=1}^N e_{\! j^k}$,
where $e_{\! j}$ is the $j$th unit vector in $\Rb^n$. Since the sample is finite, we can calculate the plug-in risk measure estimate,
\begin{equation}
\label{plug-in}
\widetilde{\sigma}^N_i(P_i,v) = \widebar{\sigma}_i(P_i^N,v),
\end{equation}
by a closed-form expression. One can verify directly from the definition that the resulting \emph{sample-based transition risk mapping}
\[
\sigma_i^N(P_i,v) = \Eb\big[ \widebar{\sigma}_i(P_i^N,v) \big],
\]
satisfies all conditions of a transition risk mapping of \S \ref{sec:projec}, if $\widebar{\sigma}_i(\cdot,\cdot)$ does. The expectation above is over all possible $N$-samples. Therefore, if we treat $\sigma_i^N(\cdot,\cdot)$ as the ``true'' risk measure that we want to estimate, the plug-in
formula \eqref{plug-in} satisfies \eqref{sigma-tilde} and Assumption \ref{a:xi}. In fact, for a broad class of measures of risk
$\widebar{\sigma}_i(P_i,v)$, we have a central limit result: $\widebar{\sigma}_i(P_i^N,v)$ is convergent to $\widebar{\sigma}_i(P_i,v)$
at the rate $1/\sqrt{N}$, and the error has an approximately normal distribution \cite{dentcheva2017statistical}. However, we do not rely on this result here, because we work with fixed $N$. In our experiments,
the sample size $N=4$ turned out to be sufficient, and even $N=2$ would work well.

\subsection{Example}

We apply the risk-averse methods of temporal differences to a version of a transportation problem discussed in \cite{Powell2006}. We have vehicles  at $M=50$ locations. At each time period $t$, a stochastic demand $D_{i\! jt}$ for transportation from location $i$ to location $j$ occurs, $i,j=1,\dots,M$, $t = 1,2,\ldots$. The demand arrays $D_t$ in different time periods
are independent.
The vehicles available at location $i$ may be used to satisfy this demand. They may also be moved empty.
The state $x_t$ of the system at time $t$ is the $M$-dimensional integer vector containing the numbers of vehicles at each location.

For simplicity, we assume that a vehicle can carry a unit demand, and the total demand at the location $i$ at time $t$ can be satisfied only if $x_{it} \geq \sum_{j =1}^M D_{i\! jt}$; otherwise, the demand may be only partially satisfied and the excess demand is lost. One can relocate the vehicles empty or loaded, and we denote the cost of moving a vehicle empty from location $i$ to location $j$ as $c^e_{ij}$. Since we stay in a cost minimization setting, we also denote the net negative profit of moving a vehicle loaded from location $i$ to location $j$ as $c^{\ell}_{ij}$. Let $u_{i\! jt}^e$ be the number of vehicles moved empty from location $i$ to location $j$ at time $t$ and $u_{i\! jt}^\ell$ be the number of vehicles that are moved loaded. For simplicity, let us refer  to the combination of $u_t^e$ and $u_t^\ell$ as $u_t$ and denote:
\begin{align*}
c^{\top}u_t = \sum_{i,j =1}^M \big( c_{ij}^e u_{i\! jt}^e + c_{ij}^\ell u_{i\! jt}^\ell\big).
\end{align*}
In this problem, the control $u_t$ is decided \emph{after} the state $x_{t}$ and the demand $D_t$ are observed.
The next state is a linear function of $x_t$ and $u_t$:
\[
x_{t+1} = x_t - A u_t,
\]
where $A$ can be written in an explicit way by counting the outgoing and incoming vehicles.

We denote by $\mathcal{U}(x_t, D_t)$ the set of decisions that can be taken at state $x_t$ under demand $D_t$. Our approach allows us to evaluate a  look-ahead policy  defined by a simple linear programming problem:
\begin{equation}
\label{u-pi}
u_t^\pi (x_t,D_t) = \argmin_{u \in \mathcal{U}(x_t, D_t)}\Big\{ c^{\top}u + \alpha {\pi}^{\top} (x_t - A u)\Big\}.
\end{equation}
Here, ${\pi}$ is the vector of approximate next-state values fully defining the policy. In our case, the immediate cost $ c^{\top}u_t $
depends on $D_t$, and thus the
risk-averse policy evaluation equation \eqref{policy-risk} has the following form:
\[
v^\pi(x) = \sigma\Big(P, c^\top u^\pi(x,D) + \alpha v^\pi \big(x - A u^\pi(x,D)\big)\Big),
\]
with $P$ denoting the distribution of the demand.
Our objective is to evaluate the policy $\pi$ and to improve it. As the size of the state space is enormous, we resort
to linear approximations of form $\eqref{v-tilde-linear}$, using the state $x$ as the feature vector:
$\widetilde{v}(x_t) =  x_t^\top r$.
The approximate risk-averse dynamic programming equation \eqref{projeq} takes on the form:
\begin{align}
r^\top x  = \sigma\Big(P, c^\top u^\pi(x,D) + \alpha r^\top \big(x - A u^\pi(x,D)\big)\Big). \label{projeq-st}
\end{align}
We omit the projection operator, because the feature space has full dimension. Thanks to that,
the multistep approximate risk-averse dynamic programming equation \eqref{projeq-L} coincides with \eqref{projeq-st},
and all risk-averse methods with $\lambda\in [0,1]$ solve the same equation.

In fact, we can combine the learning and policy improvement in one process, known as the \emph{optimistic approach}, in which we
always use the current $r_t$ as the vector $\pi$ defining the policy.

\subsection{Results}

 We tested the risk-averse and the risk-neutral TD($\lambda$) methods under the same long simulated sequence of demand vectors.
  At every time~$t$, we sampled $N=4$ instances of the demand vectors, and for each instance, we computed the best decisions by \eqref{u-pi},
  and the resulting states. Then we computed the empirical risk measure
 \eqref{plug-in}
 of the approximate value of the next state, and we used it in the observed temporal difference calculation \eqref{TD0-1}:
 \[
\widetilde{d}_t = r_t^\top x_t - \alpha \widebar{\sigma}\Big(P^N, c^\top u^{r_t}(x_t,D) + \alpha r_t^\top \big(x_t - A u^{r_t}(x_t,D)\big)\Big).
 \]
We used the mean--semideviation risk measure \cite{OR:1999} as $\widebar{\sigma}(\cdot,\cdot)$, which can be calculated in closed form for an empirical distribution $P^N$ with observed transition costs  $w^{(1)},\dots,v^{(N)}$:
 \[
 \widebar{\sigma}(P^N,v) = \mu + \beta \frac{1}{N} \sum_{j=1}^N \max(0, w^{(j)}-\mu),\quad \mu = \frac{1}{N} \sum_{j=1}^N  w^{(j)},\quad \beta\in [0,1].
 \]
 We used $\beta = 1$, $N=4$, and $\alpha=0.95$.
 In the expected value model ($\beta=0$),
 we also used $N=4$ observations per stage, and we averaged them, to make the comparison fair. The choice of $N=4$ was due to the use
 of a four-core computer, on which the $N$ transitions can be simulated and analyzed in parallel.

\begin{figure}[h]
	\centering
	\subfloat[$\lambda = 0$]{%
		\includegraphics[width=.5\textwidth]{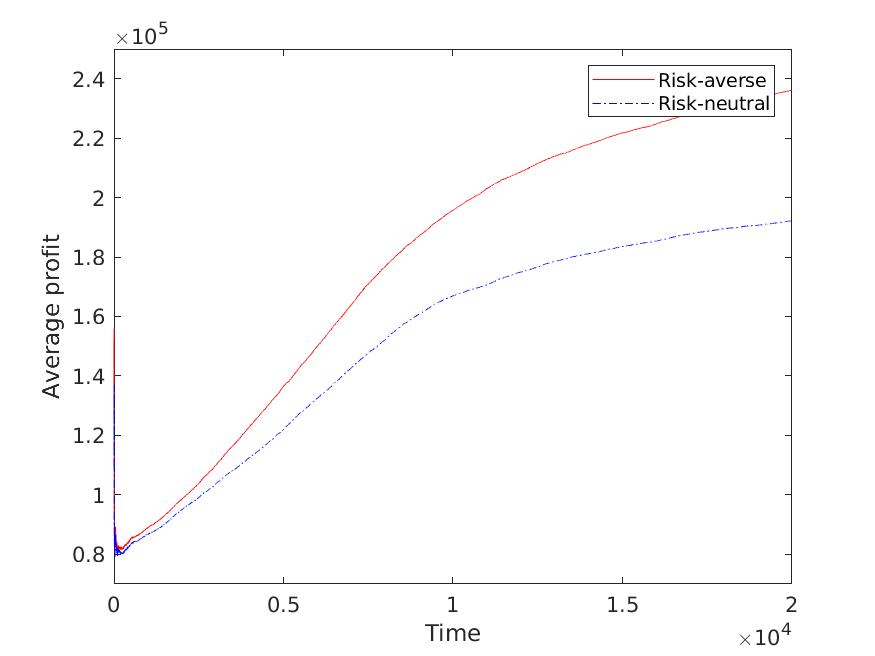}}\hfill
	\subfloat[$\lambda = 0.5$]{%
		\includegraphics[width=.5\textwidth]{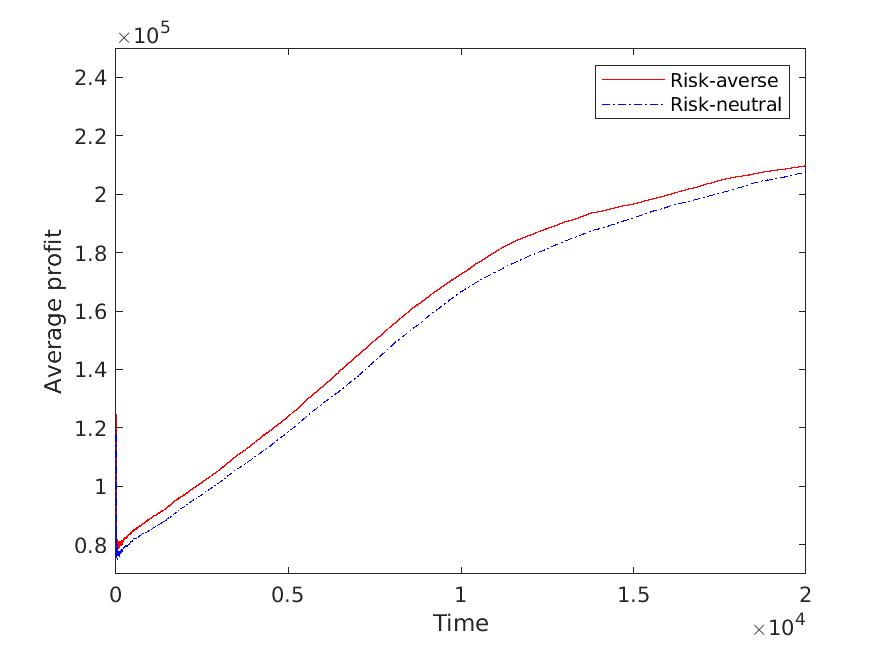}}\hfill\\
	\subfloat[$\lambda = 0.9$]{%
		\includegraphics[width=.5\textwidth]{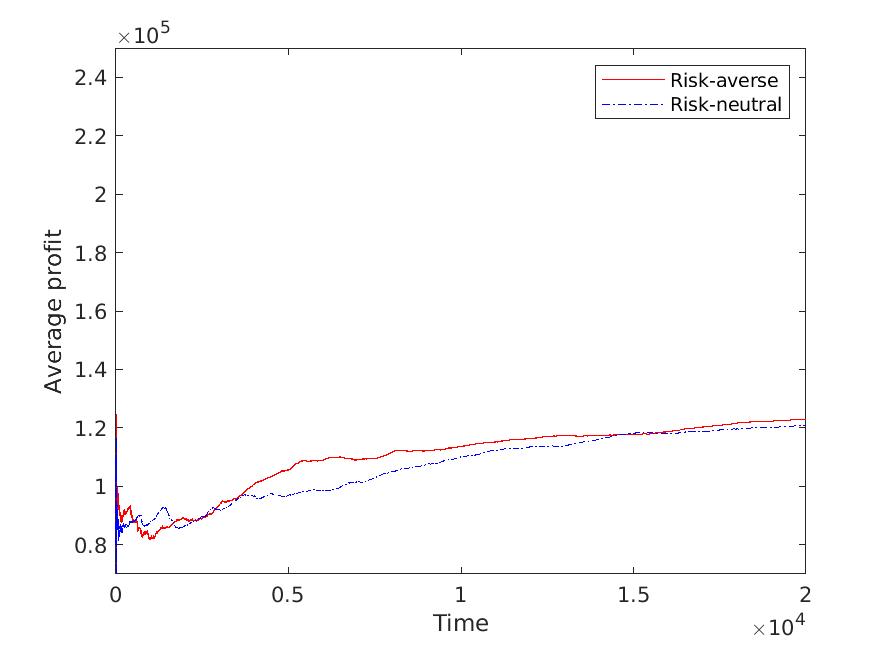}}\hfill
	\caption{Evolution of the average profit per stage.}\label{fig:avg20k} \vspace{-3ex}
\end{figure}

 We compared the performance of the risk-averse and risk-neutral TD($\lambda$) algorithms for $\lambda = 0$,   0.5, and 0.9,  in terms of average profit per stage, on a trajectory with 20,000 decision stages. The results are depicted in Figure \ref{fig:avg20k}.
We observe that the risk-averse algorithms outperform their risk-neutral counterparts in terms of the average profit in the long run. We also observe that the difference in performance is more significant when $\lambda$ is closer to zero. It would appear that
with risk-averse learning  no additional advantage is gained by using $\lambda >0$.

\begin{figure}[h]
	\centering
	\subfloat[$\lambda = 0$]{%
		\includegraphics[width=.5\textwidth]{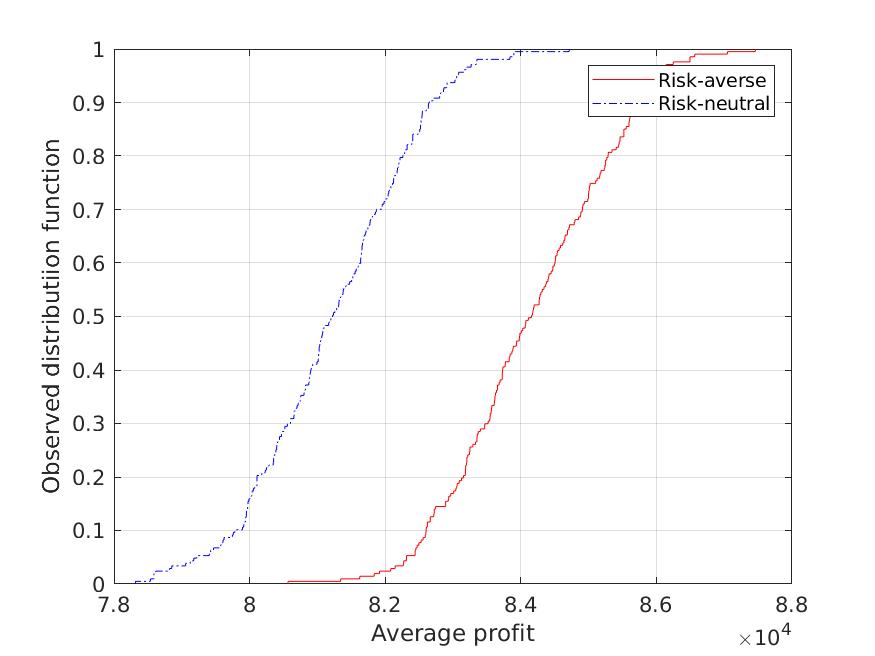}}\hfill
	\subfloat[$\lambda = 0.5$]{%
		\includegraphics[width=.5\textwidth]{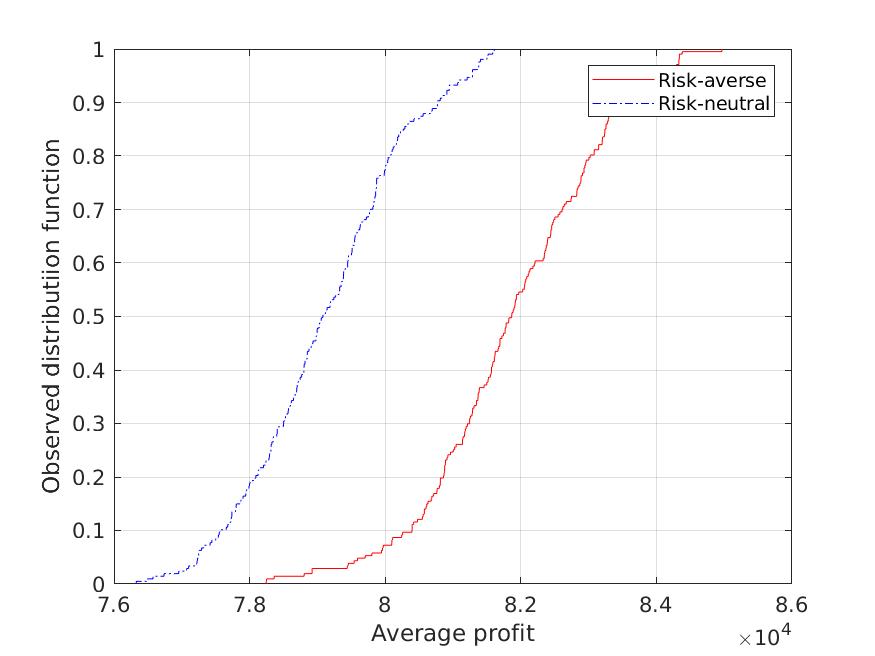}}\hfill\\
	\subfloat[$\lambda = 0.9$]{%
		\includegraphics[width=.5\textwidth]{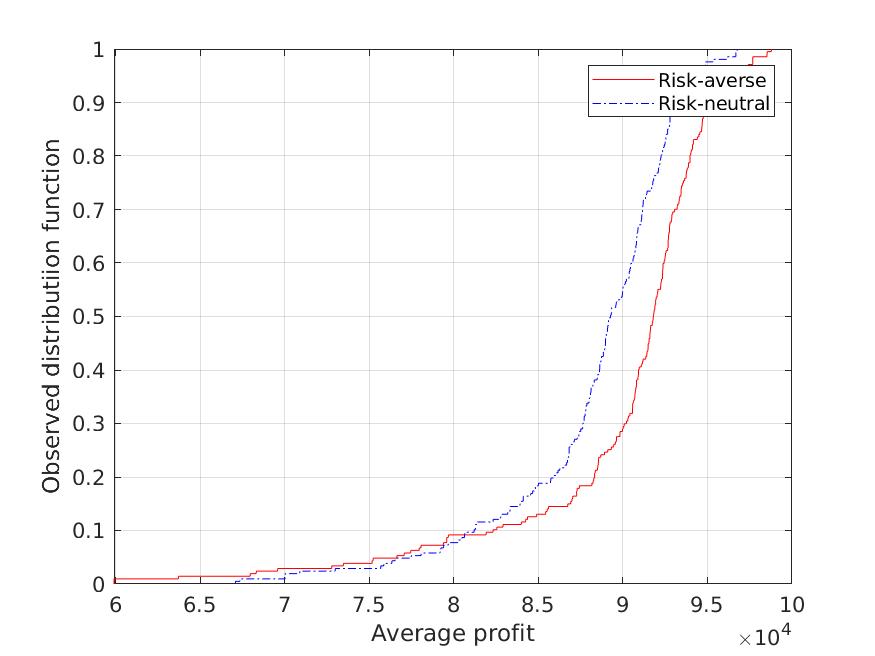}}\hfill
	\caption{Empirical distribution of the average profit at $t=200$.}\label{fig:avg200} \vspace{-3ex}
\end{figure}

In addition to these results, we used 207 distinct trajectories, each with 200 decision stages, to compare the performance of the risk-averse and risk-neutral algorithms at the early training stages in terms of profit per stage. Figure \ref{fig:avg200} shows the empirical distribution function of the profit per stage of the risk-averse and risk-neutral algorithms at $t=200$, for $\lambda = 0$, 0.5, and 0.9. The results demonstrate that in the early stages of learning ($t=200$), the average profit of the risk-averse algorithm is more likely to be higher than that of the risk-neutral algorithm, and the difference is very pronounced for lower values of $\lambda$. The first order stochastic dominance relation between empirical distributions appears to exist.

Although the risk-averse methods aim at optimizing the dynamic risk measure, rather than the expected value, they
outperform the expected value model also in expectation. This may be due to the fact that the use of risk measures makes the
method less sensitive to the imperfections of the value function approximation.

\section*{Acknowledgments}
The authors acknowledge the Office of Advanced Research Computing (http://oarc.rutgers.edu) at Rutgers, The State University of New Jersey, for providing access to the Amarel cluster and associated research computing resources that have contributed to the results reported here.
\newpage
\bibliographystyle{plain}

\end{document}